\documentclass[11pt, a4paper]{article}
\usepackage{mathrsfs, mathtools, amsthm, amssymb, thm-restate}  
\usepackage{paralist, tablists}  

\usepackage{graphicx, xcolor, tikz}  
\usetikzlibrary{calc, shapes, decorations.pathreplacing, decorations.markings}  
\usepackage{caption}  
\usepackage[labelformat=simple]{subcaption}  

\usepackage[left=25mm, top=25mm, bottom=25mm, right=25mm]{geometry}  
\usepackage[T1]{fontenc} 

\usepackage[inline]{enumitem}  
\usepackage[square, numbers, sort&compress]{natbib}  
\usepackage[
  bookmarks=true,                
  bookmarksnumbered=true,        
  bookmarksopen=true,            
  plainpages=false,              
  colorlinks=true,               
  linkcolor=blue,                
  citecolor=black,               
  anchorcolor=green,             
  urlcolor=blue,                 
  hyperindex=true                
]{hyperref} 

\newtheorem{theorem}{Theorem}[section] 
\newtheorem{lemma}[theorem]{Lemma}

\theoremstyle{definition}

\makeatletter
\def\th@plain{%
  \upshape 
}
\makeatother

\newcommand{\etal}{et~al.\ }

\usepackage{marvosym,wasysym}
\usepackage{array}
\usepackage{bm}  

\makeatletter
\renewenvironment{proof}[1][\proofname]{\par
  \pushQED{\qed}%
  \normalfont \topsep6\p@\@plus6\p@\relax
  \trivlist
  \item[\hskip\labelsep
        \bfseries
    #1\@addpunct{.}]\ignorespaces
}{%
  \popQED\endtrivlist\@endpefalse
}
\makeatother

\usepackage[capitalise]{cleveref}  
\crefname{claim}{Claim}{Claims}

\tikzset{
  on each segment/.style={
    decorate,
    decoration={
      show path construction,
      moveto code={},
      lineto code={
        \path [#1] (\tikzinputsegmentfirst) -- (\tikzinputsegmentlast);
      },
      curveto code={
        \path [#1] (\tikzinputsegmentfirst)
        .. controls (\tikzinputsegmentsupporta) and (\tikzinputsegmentsupportb) ..
        (\tikzinputsegmentlast);
      },
      closepath code={
        \path [#1] (\tikzinputsegmentfirst) -- (\tikzinputsegmentlast);
      },
    },
  },
  mid arrow/.style={postaction={decorate,decoration={
        markings,
        mark=at position .6 with {\arrow[#1]{stealth}}
      }}},
  end arrow/.style={postaction={decorate,decoration={
        markings,
        mark=at position 1 with {\arrow[#1]{stealth}}
      }}},
}

\begin{document}

\title{Decomposition of toroidal graphs without some subgraphs}
\author{Tao Wang\thanks{Center for Applied Mathematics, Henan University, Kaifeng, 475004, P. R. China} \and Xiaojing Yang\thanks{School of Mathematics and Statistics, Henan University, Kaifeng, 475004, P. R. China. {\tt Corresponding author: yangxiaojing@henu.edu.cn}. The second author was supported by National Natural Science Foundation of China (No. 12101187).}}

\date{}
\maketitle
\begin{abstract}
We consider a family of toroidal graphs, denoted by $\mathcal{T}_{i, j}$, which contain neither $i$-cycles nor $j$-cycles. A graph $G$ is $(d, h)$-decomposable if it contains a subgraph $H$ with $\Delta(H) \leq h$ such that $G - E(H)$ is a $d$-degenerate graph. For each pair $(i, j) \in \{(3, 4), (3, 6), (4, 6), (4, 7)\}$, Lu and Li proved that every graph in $\mathcal{T}_{i, j}$ is $(2, 1)$-decomposable. In this short note, we present a unified approach to prove that a common superclass of $\mathcal{T}_{i, j}$ is also $(2, 1)$-decomposable. 

\vspace*{0.2cm}\noindent {\bf Keywords}: Toroidal graph; Decomposition; Discharging

\vspace*{0.2cm}\noindent {\bf MSC2020}: 05C10; 05C15
\end{abstract}

\section{Introduction}
For a $2$-cell embedded graph $G$ on an orientable surface, the sets of the vertices, the edges, and the faces are, respectively, represented by $V(G), E(G)$, and $F(G)$. Throughout this paper, when we refer a face of a graph, we mean the graph has been embedded on a plane or a torus. For any subset $S \subseteq V(G)$, the complement of $S$ in $G$ is written as $\overline{S} = V(G) - S$. The \emph{degree} of a vertex $v$ in a graph (or the \emph{out-degree} in a digraph) is denoted by $\deg(v)$ (or $\deg^{+}(v)$ in the digraph). We use $\Delta^{+}(D)$ to denote the maximum out-degree of an orientation $D$. For a vertex $v$, if $\deg(v) = d$ (resp. $\deg(v) \leq d$, or $\deg(v) \geq d$), then we say that $v$ is a $d$-vertex (resp. $d^{-}$-vertex, or $d^{+}$-vertex). Similarly, for a face $f$ of an embedded graph, if its size is $d$ (resp. at most $d$ or at least $d$), then we call $f$ a $d$-face (resp. $d^{-}$-face or $d^{+}$-face). For a face $f$ in $F(G)$, we write $f = [v_1v_2 \dots v_m]$ to denote a cyclic order of vertices on the boundary. An $l$-face $[v_1v_2\dots v_l]$ is referred to as a $(d_1,d_2,\dots,d_l)$-face, provided $\deg(v_j)=d_j$ for all $1 \leq j \leq l$. Two faces are defined to be \emph{adjacent} if they have a common edge, and \emph{normally adjacent} if they share exactly one edge and exactly two common vertices. The symbol $\overrightarrow{uv}$ is used to indicate the orientation of the edge $uv$ from $u$ to $v$.

A \emph{$(t_{1}, t_{2}, \dots, t_{k})$-coloring} is a vertex partition $(U_{1}, U_{2}, \dots, U_{k})$ where the maximum degree of $G[U_{i}]$ is at most $t_{i}$ for all $1 \leq i \leq k$. Obviously, a $(0_{1}, 0_{2}, \dots, 0_{k})$-coloring is equivalent to a proper vertex $k$-coloring. An \emph{$h$-defective $k$-coloring} is referred to as a $k$-coloring where the maximum degree of each color class is at most $h$. An \emph{$h$-defective $L$-coloring} of a graph $G$ is defined as an $h$-defective $k$-coloring $\phi$ of $G$ satisfying $\phi(u)$ belongs to $L(u)$ for every vertex $u$ in $V(G)$. Furthermore, $G$ is said to be \emph{$h$-defective $k$-chosen} if it has an $h$-defective $L$-coloring for each $k$-list assignment $L$.

Cowen \etal \cite{MR890224} have proved that planar graphs can be $2$-defective $3$-colored. Subsequently, Eaton and Hull \cite{MR1668108}, as well as \v{S}krekovski \cite{MR1702609}, independently confirmed the $2$-defective $3$-choosability of planar graphs. In \cite{MR1668108,MR1702609}, the authors raised the problem that whether planar graphs can be $1$-defective $4$-chosen. Cushing and Kierstead \cite{MR2644426} provided a positive solution to this inquiry. 

An \emph{$l$-cycle} is a cycle with a length of $l$. Consider the family of planar graphs denoted by $\mathcal{P}_{i, j}$, which consists of planar graphs that do not contain $i$-cycles or $j$-cycles. In \cite{MR1820611}, the authors established that for every $j \in \{5, 6, 7\}$, each member in $\mathcal{P}_{4, j}$ can be $1$-defective $3$-chosen. This conclusion was also verified by Dong and Xu \cite{MR2479819} for each $l \in \{8, 9\}$.

A \emph{$(d, h)$-decomposition} of a graph $G$ is an ordered pair $(D, H)$ such that $H$ is a subgraph of $G$ with $\Delta(H) \leq h$, and $D$ is an acyclic orientation of $G - E(H)$ with $\Delta^{+}(D) \leq d$. We say a graph is \emph{$(d, h)$-decomposable} if it has a $(d, h)$-decomposition. A \emph{$d$-degenerate graph} is a graph whose every subgraph has minimum degree at most $d$. Note that a graph is $d$-degenerate if and only if it has an orientation $D$ with $\Delta^{+}(D) \leq d$. Equivalently, a graph $G$ is $(d, h)$-decomposable if it contains a subgraph $H$ with $\Delta(H) \leq h$ and the removing of $E(H)$ from $G$ results in a $d$-degenerate graph. We also say that a graph $G$ can be \emph{$(d, h)$-decomposed} if it is $(d, h)$-decomposable. Note that $k$-degenerate graphs can be $(k+1)$-colorable, $(k+1)$-choosable, $(k+1)$-AT-colorable, $(k+1)$-DP-colorable, and $(k+1)$-DP-paintable. For the definition of $k$-AT-colorable, $k$-DP-colorable and $k$-DP-paintable, we refer the reader to \cite{MR3906645,MR3686937,MR4117373}. It follows that any $(d, h)$-decomposable graph is $h$-defective $(d+1)$-choosable. Zhu \cite{MR1746462} showed that a planar graph can be $(2, 8)$-decomposed.

Cho \etal \cite{MR4472765} considered general planar graphs and obtained the following results. 
\begin{theorem}[Cho \etal \cite{MR4472765}]\label{Cho}\text{}
\begin{enumerate}
\item[(1)] Every planar graph is $(4, 1)$-decomposable.
\item[(2)] Every planar graph is $(3, 2)$-decomposable.
\item[(3)] Every planar graph is $(2, 6)$-decomposable.
\item[(4)] Not all planar graphs are $(2, 3)$-decomposable.
\end{enumerate}
\end{theorem}

It is worth noting that the results presented in \cref{Cho}(1) and (2) are optimal since a planar graph with minimum degree at least five is neither $(3, 1)$-decomposable nor $(4, 0)$-decomposable.

Lu and Zhu \cite{MR4051856} improved Lih \etal's results to the following. 
\begin{theorem}[Lu and Zhu \cite{MR4051856}]
For each $l \in \{5, 6, 7\}$, any graph belonging to $\mathcal{P}_{4, l}$ can be $(2, 1)$-decomposed. 
\end{theorem}

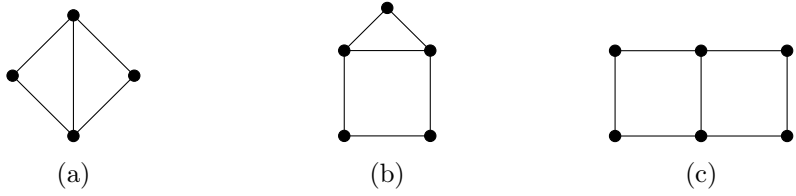
\begin{figure}
\centering
\subcaptionbox{}[0.25\linewidth]{\begin{tikzpicture}[scale=0.8]
\def\s{1}
\coordinate (E) at (\s, 0);
\coordinate (N) at (0,\s);
\coordinate (W) at (-\s,0);
\coordinate (S) at (0,-\s);
\draw (E)--(N)--(W)--(S)--cycle;
\draw (N)--(S);
\node[circle, inner sep = 1.5, fill, draw] () at (E) {};
\node[circle, inner sep = 1.5, fill, draw] () at (N) {};
\node[circle, inner sep = 1.5, fill, draw] () at (W) {};
\node[circle, inner sep = 1.5, fill, draw] () at (S) {};
\end{tikzpicture}}
\subcaptionbox{}[0.25\linewidth]{\begin{tikzpicture}[scale=0.8]
\def\s{1}
\coordinate (NE) at (45:\s);
\coordinate (NW) at (135:\s);
\coordinate (SW) at (225:\s);
\coordinate (SE) at (-45:\s);
\coordinate (H) at (90:1.414*\s);
\draw (NE)--(H)--(NW)--(SW)--(SE)--cycle;
\draw (NE)--(NW);
\node[circle, inner sep = 1.5, fill, draw] () at (NE) {};
\node[circle, inner sep = 1.5, fill, draw] () at (NW) {};
\node[circle, inner sep = 1.5, fill, draw] () at (SW) {};
\node[circle, inner sep = 1.5, fill, draw] () at (SE) {};
\node[circle, inner sep = 1.5, fill, draw] () at (H) {};
\end{tikzpicture}}
\subcaptionbox{}[0.25\linewidth]{\begin{tikzpicture}[scale=0.8]
\def\s{1.414}
\coordinate (NE) at (\s, \s);
\coordinate (N) at (0,\s);
\coordinate (NW) at (-\s,\s);
\coordinate (SW) at (-\s, 0);
\coordinate (S) at (0, 0);
\coordinate (SE) at (\s, 0);
\draw (NE)--(NW)--(SW)--(SE)--cycle;
\draw (N)--(S);
\node[circle, inner sep = 1.5, fill, draw] () at (NE) {};
\node[circle, inner sep = 1.5, fill, draw] () at (N) {};
\node[circle, inner sep = 1.5, fill, draw] () at (NW) {};
\node[circle, inner sep = 1.5, fill, draw] () at (SW) {};
\node[circle, inner sep = 1.5, fill, draw] () at (S) {};
\node[circle, inner sep = 1.5, fill, draw] () at (SE) {};
\end{tikzpicture}}
\caption{Common forbidden configurations in  \cref{LLWZ}(1) and (2).}
\label{COMMONFIGURE}
\end{figure}

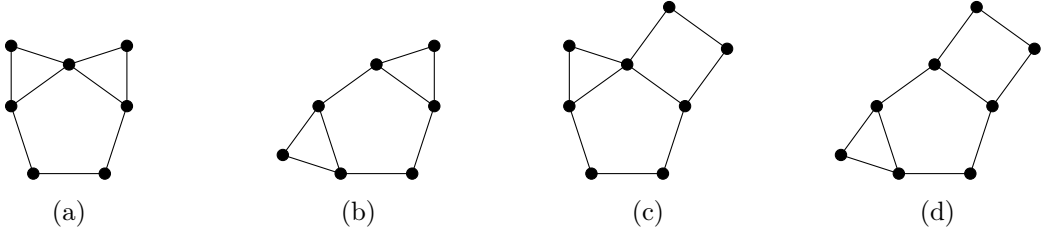
\begin{figure}
\centering
\subcaptionbox{}[0.23\linewidth]{\begin{tikzpicture}[scale=0.8]
\def\s{1}
\foreach \ang in {1, 2, 3, 4, 5}
{
\def\pointname{v\ang}
\coordinate (\pointname) at ($(\ang*360/5-54:\s)$);
}
\coordinate (NE) at ($(v1)+(v2)$);
\coordinate (NW) at ($(v2)+(v3)$);
\draw (v1)--(v2)--(v3)--(v4)--(v5)--cycle;
\draw (v1)--(NE)--(v2);
\draw (v2)--(NW)--(v3);
\foreach \ang in {1, 2, 3, 4, 5}
{
\node[circle, inner sep = 1.5, fill, draw] () at (v\ang) {};
}
\node[circle, inner sep = 1.5, fill, draw] () at (NE) {};
\node[circle, inner sep = 1.5, fill, draw] () at (NW) {};
\end{tikzpicture}}
\subcaptionbox{}[0.23\linewidth]{\begin{tikzpicture}[scale=0.8]
\def\s{1}
\foreach \ang in {1, 2, 3, 4, 5}
{
\def\pointname{v\ang}
\coordinate (\pointname) at ($(\ang*360/5-54:\s)$);
}
\coordinate (NE) at ($(v1)+(v2)$);
\coordinate (SW) at ($(v3)+(v4)$);
\draw (v1)--(v2)--(v3)--(v4)--(v5)--cycle;
\draw (v1)--(NE)--(v2);
\draw (v3)--(SW)--(v4);
\foreach \ang in {1, 2, 3, 4, 5}
{
\node[circle, inner sep = 1.5, fill, draw] () at (v\ang) {};
}
\node[circle, inner sep = 1.5, fill, draw] () at (NE) {};
\node[circle, inner sep = 1.5, fill, draw] () at (SW) {};
\end{tikzpicture}}
\subcaptionbox{}[0.23\linewidth]{\begin{tikzpicture}[scale=0.8]
\def\s{1}
\foreach \ang in {1, 2, 3, 4, 5}
{
\def\pointname{v\ang}
\coordinate (\pointname) at ($(\ang*360/5-54:\s)$);
}
\coordinate (NE2) at ($(v2)!1!90:(v1)$);
\coordinate (NE1) at ($(v1)!1!-90:(v2)$);
\coordinate (NW) at ($(v2)+(v3)$);
\draw (v1)--(v2)--(v3)--(v4)--(v5)--cycle;
\draw (v2)--(NW)--(v3);
\draw (v2)--(NE2)--(NE1)--(v1);
\foreach \ang in {1, 2, 3, 4, 5}
{
\node[circle, inner sep = 1.5, fill, draw] () at (v\ang) {};
}
\node[circle, inner sep = 1.5, fill, draw] () at (NE1) {};
\node[circle, inner sep = 1.5, fill, draw] () at (NE2) {};
\node[circle, inner sep = 1.5, fill, draw] () at (NW) {};
\end{tikzpicture}}
\subcaptionbox{}[0.23\linewidth]{\begin{tikzpicture}[scale=0.8]
\def\s{1}
\foreach \ang in {1, 2, 3, 4, 5}
{
\def\pointname{v\ang}
\coordinate (\pointname) at ($(\ang*360/5-54:\s)$);
}
\coordinate (NE2) at ($(v2)!1!90:(v1)$);
\coordinate (NE1) at ($(v1)!1!-90:(v2)$);
\coordinate (SW) at ($(v3)+(v4)$);
\draw (v1)--(v2)--(v3)--(v4)--(v5)--cycle;
\draw (v3)--(SW)--(v4);
\draw (v2)--(NE2)--(NE1)--(v1);
\foreach \ang in {1, 2, 3, 4, 5}
{
\node[circle, inner sep = 1.5, fill, draw] () at (v\ang) {};
}
\node[circle, inner sep = 1.5, fill, draw] () at (NE1) {};
\node[circle, inner sep = 1.5, fill, draw] () at (NE2) {};
\node[circle, inner sep = 1.5, fill, draw] () at (SW) {};
\end{tikzpicture}}
\caption{Certain forbidden configurations in \cref{LLWZ}(1).}
\label{FIGURE-AT567}
\end{figure}

Li \etal \cite{MR4557782} further extended the above result.

\begin{theorem}[Li \etal \cite{MR4557782}]\label{LLWZ}
If $G$ is a plane graph and any one of the three conditions below holds, then $G$ is $(2, 1)$-decomposable:
\begin{enumerate}[label = (\alph*)]
    \item $G$ does not contain any subgraph as shown in \cref{COMMONFIGURE,FIGURE-AT567}.
    \item $G$ does not contain any subgraph as shown in \cref{COMMONFIGURE,FIGURE-AT48}.
    \item $G$ is included in $\mathcal{P}_{4, 9}$.
\end{enumerate}
\end{theorem}

\begin{figure}
\centering
\subcaptionbox{}[0.19\linewidth]{\begin{tikzpicture}[scale=0.8]
\def\s{1.414}
\coordinate (A) at (\s, \s);
\coordinate (B) at (0,\s);
\coordinate (C) at (-\s,\s);
\coordinate (D) at (-\s, 0);
\coordinate (E) at (0, 0);
\coordinate (F) at (\s, 0);
\coordinate (H1) at (0.5*\s, 1.5*\s);
\coordinate (H2) at (-0.5*\s, 1.5*\s);
\draw (A)--(H1)--(B)--(H2)--(C)--(D)--(E)--(E)--(F)--cycle;
\draw (B)--(E);
\node[circle, inner sep = 1.5, fill, draw] () at (A) {};
\node[circle, inner sep = 1.5, fill, draw] () at (B) {};
\node[circle, inner sep = 1.5, fill, draw] () at (C) {};
\node[circle, inner sep = 1.5, fill, draw] () at (D) {};
\node[circle, inner sep = 1.5, fill, draw] () at (E) {};
\node[circle, inner sep = 1.5, fill, draw] () at (F) {};
\node[circle, inner sep = 1.5, fill, draw] () at (H1) {};
\node[circle, inner sep = 1.5, fill, draw] () at (H2) {};
\end{tikzpicture}}
\subcaptionbox{}[0.19\linewidth]{\begin{tikzpicture}[scale=0.8]
\def\s{1.414}
\coordinate (NE) at (0.5*\s, 0.5*\s);
\coordinate (NW) at (-0.5*\s, 0.5*\s);
\coordinate (SW) at (-0.5*\s, -0.5*\s);
\coordinate (SE) at (0.5*\s, -0.5*\s);
\coordinate (H1) at (0, \s);
\coordinate (H2) at (0, 0.6*\s);
\draw (NE)--(H1)--(NW)--(SW)--(SE)--cycle;
\draw (NE)--(H2)--(NW);
\node[circle, inner sep = 1.5, fill, draw] () at (NE) {};
\node[circle, inner sep = 1.5, fill, draw] () at (NW) {};
\node[circle, inner sep = 1.5, fill, draw] () at (SE) {};
\node[circle, inner sep = 1.5, fill, draw] () at (SW) {};
\node[circle, inner sep = 1.5, fill, draw] () at (H1) {};
\node[circle, inner sep = 1.5, fill, draw] () at (H2) {};
\end{tikzpicture}}
\subcaptionbox{}[0.19\linewidth]{\begin{tikzpicture}[scale=0.8]
\def\s{1}
\foreach \ang in {1, 2, 3, 4, 5, 6, 7}
{
\def\pointname{v\ang}
\coordinate (\pointname) at ($(\ang*360/7-90:\s)$);
}
\coordinate (W) at ($(v3)!1!-60:(v4)$);
\draw (v1)--(v2)--(v3)--(v4)--(v5)--(v6)--(v7)--cycle;
\draw (v3)--(W)--(v4);
\foreach \ang in {1, 2, 3, 4, 5, 6, 7}
{
\node[circle, inner sep = 1.5, fill, draw] () at (v\ang) {};
}
\node[circle, inner sep = 1.5, fill, draw] () at (W) {};
\end{tikzpicture}}
\subcaptionbox{}[0.19\linewidth]{\begin{tikzpicture}[scale=0.8]
\def\s{1}
\foreach \ang in {1, 2, 3, 4, 5, 6}
{
\def\pointname{v\ang}
\coordinate (\pointname) at ($(\ang*360/6:\s)$);
}
\coordinate (N) at ($(v1)+(v2)$);
\coordinate (SW1) at ($(v3)!1!-90:(v4)$);
\coordinate (SW2) at ($(v4)!1!90:(v3)$);
\draw (v1)--(v2)--(v3)--(v4)--(v5)--(v6)--cycle;
\draw (v1)--(N)--(v2);
\draw (v3)--(SW1)--(SW2)--(v4);
\foreach \ang in {1, 2, 3, 4, 5, 6}
{
\node[circle, inner sep = 1.5, fill, draw] () at (v\ang) {};
}
\node[circle, inner sep = 1.5, fill, draw] () at (N) {};
\node[circle, inner sep = 1.5, fill, draw] () at (SW1) {};
\node[circle, inner sep = 1.5, fill, draw] () at (SW2) {};
\end{tikzpicture}}
\subcaptionbox{}[0.19\linewidth]{\begin{tikzpicture}[scale=0.8]
\def\s{1}
\foreach \ang in {1, 2, 3, 4, 5, 6}
{
\def\pointname{v\ang}
\coordinate (\pointname) at ($(\ang*360/6:\s)$);
}
\coordinate (W) at ($(v1)+(v2)$);
\coordinate (S1) at ($(v4)!1!-90:(v5)$);
\coordinate (S2) at ($(v5)!1!90:(v4)$);
\draw (v1)--(v2)--(v3)--(v4)--(v5)--(v6)--cycle;
\draw (v1)--(W)--(v2);
\draw (v4)--(S1)--(S2)--(v5);
\foreach \ang in {1, 2, 3, 4, 5, 6}
{
\node[circle, inner sep = 1.5, fill, draw] () at (v\ang) {};
}
\node[circle, inner sep = 1.5, fill, draw] () at (W) {};
\node[circle, inner sep = 1.5, fill, draw] () at (S1) {};
\node[circle, inner sep = 1.5, fill, draw] () at (S2) {};
\end{tikzpicture}}
\caption{Certain forbidden configurations in \cref{LLWZ}(2).}
\label{FIGURE-AT48}
\end{figure}

Recently, some scholars focused on studying toroidal graphs. 
\begin{theorem}[Lu and Li \cite{Lu2023}]
Every toroidal graph without adjacent triangles is $(3, 1)$-decomposable.
\end{theorem}

\begin{theorem}[Tian \etal \cite{MR4593925}]
A toroidal graph can be $(3, 1)$-decomposed if one of the following conditions holds:
\begin{enumerate}
\item[(1)] every $5$-cycle is an induced cycle;
\item[(2)] every $6$-cycle is an induced cycle. 
\end{enumerate}
\end{theorem}

Consider the family of toroidal graphs denoted by $\mathcal{T}_{i, j}$, which consists of toroidal graphs that do not contain $i$-cycles or $j$-cycles. Lu and Li \cite{Lu2023} focused on this family of toroidal graphs. 

\begin{theorem}[Lu and Li \cite{Lu2023}]\label{Lu}
For each pair $(i, j) \in \{(3, 4), (3, 6), (4, 6), (4, 7)\}$, every graph in $\mathcal{T}_{i, j}$ can be $(2, 1)$-decomposed. 
\end{theorem}

\begin{figure}
\centering
\subcaptionbox{\label{CommonA}}[0.15\linewidth]{\begin{tikzpicture}
\def\s{1}
\coordinate (E) at (\s, 0);
\coordinate (N) at (0,\s);
\coordinate (W) at (-\s,0);
\coordinate (S) at (0,-\s);
\draw (E)--(N)--(W)--(S)--cycle;
\draw (N)--(S);
\node[circle, inner sep =1.5, fill, draw] () at (E) {};
\node[circle, inner sep =1.5, fill, draw] () at (N) {};
\node[circle, inner sep =1.5, fill, draw] () at (W) {};
\node[circle, inner sep =1.5, fill, draw] () at (S) {};
\end{tikzpicture}}
\subcaptionbox{\label{CommonB}}[0.15\linewidth]{\begin{tikzpicture}
\def\s{1}
\coordinate (NE) at (45:\s);
\coordinate (NW) at (135:\s);
\coordinate (SW) at (225:\s);
\coordinate (SE) at (-45:\s);
\coordinate (H) at ($(NW)!1!60:(NE)$);
\draw (NE)--(H)--(NW)--(SW)--(SE)--cycle;
\draw (NE)--(NW);
\node[circle, inner sep =1.5, fill, draw] () at (NE) {};
\node[circle, inner sep =1.5, fill, draw] () at (NW) {};
\node[circle, inner sep =1.5, fill, draw] () at (SW) {};
\node[circle, inner sep =1.5, fill, draw] () at (SE) {};
\node[circle, inner sep =1.5, fill, draw] () at (H) {};
\end{tikzpicture}}
\subcaptionbox{\label{CommonC}}[0.15\linewidth]{\begin{tikzpicture}
\def\s{1.4}
\coordinate (NE) at (0.5*\s, \s);
\coordinate (NW) at (-0.5*\s,\s);
\coordinate (SW) at (-0.5*\s, -\s);
\coordinate (SE) at (0.5*\s, -\s);
\coordinate (E) at (0.5*\s, 0);
\coordinate (W) at (-0.5*\s, 0);
\draw (NE)--(NW)--(SW)--(SE)--cycle;
\draw (E)--(W);
\node[circle, inner sep =1.5, fill, draw] () at (NE) {};
\node[circle, inner sep =1.5, fill, draw] () at (NW) {};
\node[circle, inner sep =1.5, fill, draw] () at (SW) {};
\node[circle, inner sep =1.5, fill, draw] () at (SE) {};
\node[circle, inner sep =1.5, fill, draw] () at (E) {};
\node[circle, inner sep =1.5, fill, draw] () at (W) {};
\end{tikzpicture}}
\subcaptionbox{\label{AT48D}}[0.15\linewidth]{\begin{tikzpicture}
\def\s{1}
\foreach \ang in {1, 2, 3, 4, 5, 6}
{
\def\pointname{v\ang}
\coordinate (\pointname) at ($(\ang*360/6:\s)$);
\node[circle, inner sep =1.5, fill, draw] () at (v\ang) {};}
\coordinate (N) at ($(v1)+(v2)$);
\draw (v1)--(v2)--(v3)--(v4)--(v5)--(v6)--cycle;
\draw (v1)--(N)--(v2);
\node[circle, inner sep =1.5, fill, draw] () at (N) {};
\end{tikzpicture}}
\subcaptionbox{\label{AT345C}}[0.15\linewidth]{\begin{tikzpicture}
\def\s{1}
\foreach \ang in {1, 2, 3, 4, 5}
{
\def\pointname{v\ang}
\coordinate (\pointname) at ($(\ang*360/5-54:\s)$);
\node[circle, inner sep =1.5, fill, draw] () at (v\ang) {};}
\coordinate (NE2) at ($(v2)!1!90:(v1)$);
\coordinate (NE1) at ($(v1)!1!-90:(v2)$);
\coordinate (NW) at ($(v2)+(v3)$);
\draw (v1)--(v2)--(v3)--(v4)--(v5)--cycle;
\draw (v2)--(NW)--(v3);
\draw (v2)--(NE2)--(NE1)--(v1);
\node[circle, inner sep =1.5, fill, draw] () at (NE1) {};
\node[circle, inner sep =1.5, fill, draw] () at (NE2) {};
\node[circle, inner sep =1.5, fill, draw] () at (NW) {};
\end{tikzpicture}}
\subcaptionbox{\label{AT345D}}[0.15\linewidth]{\begin{tikzpicture}
\def\s{1}
\foreach \ang in {1, 2, 3, 4, 5}
{
\def\pointname{v\ang}
\coordinate (\pointname) at ($(\ang*360/5-54:\s)$);
\node[circle, inner sep =1.5, fill, draw] () at (v\ang) {};}
\coordinate (NE2) at ($(v2)!1!90:(v1)$);
\coordinate (NE1) at ($(v1)!1!-90:(v2)$);
\coordinate (SW) at ($(v3)+(v4)$);
\draw (v1)--(v2)--(v3)--(v4)--(v5)--cycle;
\draw (v3)--(SW)--(v4);
\draw (v2)--(NE2)--(NE1)--(v1);
\node[circle, inner sep =1.5, fill, draw] () at (NE1) {};
\node[circle, inner sep =1.5, fill, draw] () at (NE2) {};
\node[circle, inner sep =1.5, fill, draw] () at (SW) {};
\end{tikzpicture}}
\caption{Certain forbidden configurations in \cref{J}.}
\label{FIGURE-Z}
\end{figure}

In this short note, we present a result that enhances the findings of the above theorem. 
\begin{theorem}\label{J}
If $G$ is a toroidal graph without subgraphs isomorphic to any configuration as depicted in \cref{FIGURE-Z}, then $G$ is $(2, 1)$-decomposable. 
\end{theorem}

\section{Proof}
This section is devoted to providing a proof for \cref{J}. 

\begin{figure}
\centering
\subcaptionbox{\label{53Fig3}}[0.23\linewidth]{\begin{tikzpicture}[line width = 1pt]
\def\s{1.4}
\coordinate (O) at (0, 0);
\coordinate (v1) at (60:\s);
\coordinate (v2) at (120:\s);
\coordinate (v3) at (270:\s);
\draw
(v2)--(O)
(O)--(v1)
(v1)--($(v1)+(0:0.5*\s)$)
(v1)--($(v1)+(60:0.5*\s)$)
(v2)--($(v2)+(150:0.5*\s)$)
(O)--($(O)+(0:0.5*\s)$)
(v3)--($(v3)+(225:0.5*\s)$)
(v3)--($(v3)+(-45:0.5*\s)$)
;
\draw
(v1)--(v2)
(O)--(v3);
\node[rectangle, inner sep = 2.5, fill, draw] () at (v1) {};
\node[rectangle, inner sep = 2.5, fill, draw] () at (O) {};
\node[regular polygon, regular polygon sides=3, inner sep = 1, fill, draw] () at (v2) {};
\node[regular polygon, regular polygon sides=3, inner sep = 1, fill, draw] () at (v3) {};
\end{tikzpicture}}
\subcaptionbox{\label{53Fig1}}[0.23\linewidth]{\begin{tikzpicture}[line width = 1pt]
\def\s{1}
\foreach \ang in {1, 2, 3, 4, 5}
{
\def\pointname{v\ang}
\coordinate (\pointname) at ($(\ang*360/5-18:\s)$);}
\coordinate (H) at ($(v2)!1!60:(v1)$);
\draw
(v3)--(v4)
(v5)--(v1)--(v2)--(H)
(v1)--($(v1)+(30:0.5*\s)$)
(v2)--($(v2)+(150:0.5*\s)$)
(H)--($(H)+(60:0.5*\s)$)
(H)--($(H)+(120:0.5*\s)$)
(v3)--($(v3)+(210:0.5*\s)$)
(v5)--($(v5)+(-30:0.5*\s)$)
(v4)--($(v4)+(-45:0.5*\s)$)
(v4)--($(v4)+(225:0.5*\s)$)
;
\draw
(v1)--(H)
(v2)--(v3)
(v4)--(v5);
\node[rectangle, inner sep = 2.5, fill, draw] () at (v1) {};
\node[rectangle, inner sep = 2.5, fill, draw] () at (v2) {};
\node[rectangle, inner sep = 2.5, fill, draw] () at (H) {};
\node[rectangle, inner sep = 2.5, fill, draw] () at (v4) {};
\node[regular polygon, regular polygon sides=3, inner sep = 1, fill, draw] () at (v3) {};
\node[regular polygon, regular polygon sides=3, inner sep = 1, fill, draw] () at (v5) {};
\end{tikzpicture}}
\subcaptionbox{\label{53Fig2}}[0.23\linewidth]{\begin{tikzpicture}[line width = 1pt]
\def\s{1}
\foreach \ang in {1, 2, 3, 4, 5}
{
\def\pointname{v\ang}
\coordinate (\pointname) at ($(\ang*360/5-18:\s)$);}
\coordinate (H) at ($(v2)!1!60:(v1)$);
\draw
(v2)--(H)
(v2)--(v1)--(v5)
(v4)--(v3)
(v1)--($(v1)+(30:0.5*\s)$)
(H)--($(H)+(60:0.5*\s)$)
(H)--($(H)+(120:0.5*\s)$)
(v3)--($(v3)+(150:0.5*\s)$)
(v3)--($(v3)+(210:0.5*\s)$)
(v5)--($(v5)+(-30:0.5*\s)$)
(v5)--($(v5)+(30:0.5*\s)$)
(v4)--($(v4)+(-90:0.5*\s)$)
;
\draw
(v1)--(H)
(v2)--(v3)
(v4)--(v5);
\node[rectangle, inner sep = 2.5, fill, draw] () at (v1) {};
\node[regular polygon, regular polygon sides=3, inner sep = 1, fill, draw] () at (v2) {};
\node[rectangle, inner sep = 2.5, fill, draw] () at (H) {};
\node[regular polygon, regular polygon sides=3, inner sep = 1, fill, draw] () at (v4) {};
\node[rectangle, inner sep = 2.5, fill, draw] () at (v3) {};
\node[rectangle, inner sep = 2.5, fill, draw] () at (v5) {};
\end{tikzpicture}}
\caption{Reducible configurations.}
\label{REDUCIBLEFIGURE}
\end{figure}
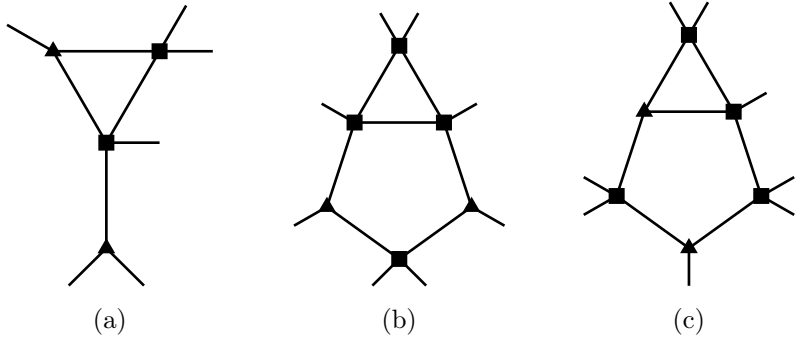

\begin{figure}
\centering
\subcaptionbox{}[0.23\linewidth]{\begin{tikzpicture}[line width = 1pt]
\def\s{1.4}
\coordinate (O) at (0, 0);
\coordinate (v1) at (60:\s);
\coordinate (v2) at (120:\s);
\coordinate (v3) at (270:\s);
\path [draw=red, postaction={on each segment={mid arrow=red}}]
(v2)--(O)
(O)--(v1);
\path [draw=red, postaction={on each segment={end arrow=red}}]
(v1)--($(v1)+(0:0.5*\s)$)
(v1)--($(v1)+(60:0.5*\s)$)
(v2)--($(v2)+(150:0.5*\s)$)
(O)--($(O)+(0:0.5*\s)$)
(v3)--($(v3)+(225:0.5*\s)$)
(v3)--($(v3)+(-45:0.5*\s)$)
;
\draw[line width = 1.5pt]
(v1)--(v2)
(O)--(v3);
\node[rectangle, inner sep = 2.5, fill, draw] () at (v1) {};
\node[rectangle, inner sep = 2.5, fill, draw] () at (O) {};
\node[regular polygon, regular polygon sides=3, inner sep = 1, fill, draw] () at (v2) {};
\node[regular polygon, regular polygon sides=3, inner sep = 1, fill, draw] () at (v3) {};
\end{tikzpicture}}
\subcaptionbox{}[0.23\linewidth]{\begin{tikzpicture}[line width = 1pt]
\def\s{1}
\foreach \ang in {1, 2, 3, 4, 5}
{
\def\pointname{v\ang}
\coordinate (\pointname) at ($(\ang*360/5-18:\s)$);}
\coordinate (H) at ($(v2)!1!60:(v1)$);
\path [draw=red, postaction={on each segment={mid arrow=red}}]
(v3)--(v4)
(v5)--(v1)--(v2)--(H);
\path [draw=red, postaction={on each segment={end arrow=red}}]
(v1)--($(v1)+(30:0.5*\s)$)
(v2)--($(v2)+(150:0.5*\s)$)
(H)--($(H)+(60:0.5*\s)$)
(H)--($(H)+(120:0.5*\s)$)
(v3)--($(v3)+(210:0.5*\s)$)
(v5)--($(v5)+(-30:0.5*\s)$)
(v4)--($(v4)+(-45:0.5*\s)$)
(v4)--($(v4)+(225:0.5*\s)$)
;
\draw[line width = 1.5pt]
(v1)--(H)
(v2)--(v3)
(v4)--(v5);
\node[rectangle, inner sep = 2.5, fill, draw] () at (v1) {};
\node[rectangle, inner sep = 2.5, fill, draw] () at (v2) {};
\node[rectangle, inner sep = 2.5, fill, draw] () at (H) {};
\node[rectangle, inner sep = 2.5, fill, draw] () at (v4) {};
\node[regular polygon, regular polygon sides=3, inner sep = 1, fill, draw] () at (v3) {};
\node[regular polygon, regular polygon sides=3, inner sep = 1, fill, draw] () at (v5) {};
\end{tikzpicture}}
\subcaptionbox{}[0.23\linewidth]{\begin{tikzpicture}[line width = 1pt]
\def\s{1}
\foreach \ang in {1, 2, 3, 4, 5}
{
\def\pointname{v\ang}
\coordinate (\pointname) at ($(\ang*360/5-18:\s)$);}
\coordinate (H) at ($(v2)!1!60:(v1)$);
\path [draw=red, postaction={on each segment={mid arrow=red}}]
(v2)--(H)
(v2)--(v1)--(v5)
(v4)--(v3);
\path [draw=red, postaction={on each segment={end arrow=red}}]
(v1)--($(v1)+(30:0.5*\s)$)
(H)--($(H)+(60:0.5*\s)$)
(H)--($(H)+(120:0.5*\s)$)
(v3)--($(v3)+(150:0.5*\s)$)
(v3)--($(v3)+(210:0.5*\s)$)
(v5)--($(v5)+(-30:0.5*\s)$)
(v5)--($(v5)+(30:0.5*\s)$)
(v4)--($(v4)+(-90:0.5*\s)$)
;
\draw[line width = 1.5pt]
(v1)--(H)
(v2)--(v3)
(v4)--(v5);
\node[rectangle, inner sep = 2.5, fill, draw] () at (v1) {};
\node[regular polygon, regular polygon sides=3, inner sep = 1, fill, draw] () at (v2) {};
\node[rectangle, inner sep = 2.5, fill, draw] () at (H) {};
\node[regular polygon, regular polygon sides=3, inner sep = 1, fill, draw] () at (v4) {};
\node[rectangle, inner sep = 2.5, fill, draw] () at (v3) {};
\node[rectangle, inner sep = 2.5, fill, draw] () at (v5) {};
\end{tikzpicture}}
\caption{Bold edges represent matching, and red edges are oriented.}
\label{DM}
\end{figure}
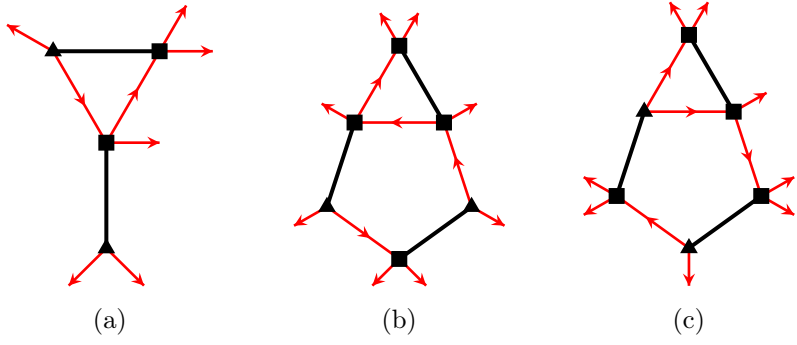

A vertex of degree $3$ is referred to as \emph{light} where it is on a $(3, 4, 4, 4)$-face and a pair of $(3, 4, 3, 4, 4)$-faces. A vertex of degree $3$ is called \emph{minor} where it is on a $4$-face. In \cref{REDUCIBLEFIGURE,XX,DM}, a triangle represents a 3-vertex, whereas a quadrilateral represents a 4-vertex.

Assume that $G$ is a counterexample to \cref{J} with minimum number of vertices, and it has been embedded on a torus. Under this assumption, $G$ has the following properties.
\begin{lemma}\label{LEM}
\mbox{}
\begin{enumerate}[label= (\roman*)]
\item\label{delta} $\delta(G)$ is at least $3$.
\item\label{a3} There are no adjacent $3$-vertices. 
\item\label{chorded} Every $5^{-}$-cycle has no chord. 
\item\label{adj4-} There are no adjacent $4^{-}$-faces. 
\item\label{5ADJACENT3} If a $5$-face is adjacent to a $4^{-}$-face, then they are normally adjacent. Moreover, every $5$-face is adjacent to at most one $3$-face.
\item\label{53ADJACENT4} If a $5$-face $f$ is adjacent to a $3$-face, then $f$ is not adjacent to any $4$-face.
\item\label{6ADJACENT3} There are no $6$-faces adjacent to $3$-faces.
\item\label{3444FACE} There are no $(3, 4, 3, 4)$-faces.
\item\label{LIGHT3VERTEX} There are no light $3$-vertices.
\item\label{ReducibleConfiguration} There are no subgraphs isomorphic to any configuration in \cref{REDUCIBLEFIGURE}.
\end{enumerate}
\end{lemma}

\begin{proof}
\ref{delta}. Suppose that a $2^{-}$-vertex $v$ exists. By minimality assumption, $G-v$ can be $(2, 1)$-decomposed into $(D', M)$. Let $D$ be defined as $D = D' \cup \{\overrightarrow{vu} \mid u\in N_{G}(v)\}$. Therefore, $G$ can be decomposed into $(D, M)$, which contradicts our assumption. 

\ref{a3}. Suppose $uv \in E(G)$ and $\deg(u) = \deg(v) = 3$. Again, by minimality assumption, $G - \{u, v\}$ admits a $(2, 1)$-decomposition $(D', M')$. Let $M$ be the union of $M'$ and $\{uv\}$, and let $D$ be obtained from $D'$ by orienting the four edges incident with $u$ and $v$ such that $\deg_{D}^{+}(u) = \deg_{D}^{+}(v) = 0$. Hence, we acquire the needed $(2, 1)$-decomposition of $G$, contradicting our initial assumption. 

\ref{chorded}. Suppose that there exists a $5^{-}$-cycle with a chord. This implies the existence of a subgraph isomorphic to \cref{CommonA} or \cref{CommonB}, which contradicts our hypothesis.

\ref{adj4-}. Since $\delta(G) \geq 3$, any two adjacent $3$-faces must be normally adjacent. Due to the absence of \cref{CommonA}, there are no adjacent $3$-faces. Assume a $4$-face $f = [uvxy]$ is adjacent to a $4^{-}$-face $g = [x'uvy'\dots]$. Again, by the absence of \cref{CommonA}, we have $y' \neq y$ and $x' \neq x$. The condition $\delta(G) \geq 3$ implies that $\{x', y'\} \cap \{x, y\} = \emptyset$. Consequently, $f$, $g$ must be normally adjacent. However, this contradicts the missing of \cref{CommonB,CommonC}.

\ref{5ADJACENT3}. Assume a $5$-face $f = [uvxyz]$ is adjacent to a $4^{-}$-face $g= [auvb \dots]$. Note that $a \neq v$ and it is possible that $a = b$. By \ref{chorded}, we have $a \notin \{x, y\}$. Since $\delta(G) \geq 3$, we also have $a \neq z$. Similarly, $b$ cannot be in $\{x, y, z, u, v\}$. Hence, $f$ and $g$ must be normally adjacent. The absence of \cref{CommonA,CommonB} implies that every vertex outside $\{u, v, x, y, z\}$ can connect with at most two neighbors in this set. If $f$ is adjacent to two $3$-faces $h_{1}$ and $h_{2}$, then the other two vertices on $h_{1}$ and $h_{2}$, not on $f$, must be distinct, and $G[V(f) \cup V(h_{1}) \cup V(h_{2})]$ contains a subgraph isomorphic to \cref{AT48D}, which is a contradiction. Therefore, every $5$-face can have at most one adjacent $3$-face. 

\ref{53ADJACENT4}. Suppose $f = [uvxyz]$ and $g = [uvw]$ have a common $uv$. By \ref{5ADJACENT3}, $g$ and $f$ are normally adjacent. Assume a $4$-face $h$ is adjacent to $f$. It is observed that $h$ and $f$ are normally adjacent. Depending on the positions of the common edges, we consider two cases. (1) Suppose $h = [uzpq]$. If $w = p$, then there exists a subgraph as shown in \cref{CommonA}, which is a contradiction. If $w = q$, then a subgraph as shown in \cref{CommonB} appears, which is a contradiction. Therefore, $w \notin \{p, q\}$, but there exists a subgraph as shown in \cref{AT345C}, which is a contradiction. (2) Suppose $h = [zypq]$. If $w \in \{p, q\}$, then there exists a subgraph isomorphic to \cref{CommonB}, which is a contradiction. Therefore, $w \notin \{p, q\}$, but there exists a subgraph isomorphic to \cref{AT345D}, which is a contradiction.

\begin{figure}
\centering
\subcaptionbox{}[0.45\linewidth]{\begin{tikzpicture}[line width = 1pt]
\def\s{1}
\coordinate (O) at (0, 0);
\coordinate (LC) at (144:\s);
\coordinate (H) at ($(LC)!1!72:(O)$);
\coordinate (L1) at ($(LC)!1!72:(H)$);
\coordinate (L2) at ($(LC)!1!72:(L1)$);
\coordinate (L3) at ($(LC)!1!72:(L2)$);
\coordinate (RC) at (36:\s);
\coordinate (R1) at ($(RC)!1!-72:(H)$);
\coordinate (R2) at ($(RC)!1!-72:(R1)$);
\coordinate (R3) at ($(RC)!1!-72:(R2)$);
\coordinate (B) at ($(R3)!1!36:(O)$);

\path [draw=red, postaction={on each segment={mid arrow=red}}]
(H)--(O)
(O)node[above right]{$a_{5}$}--(L3)
(O)--(R3)--(B)
(L1)--(H)
(L3)--(L2)
(R1)--(R2);
\path [draw=red, postaction={on each segment={end arrow=red}}]
(L1)node[left]{$a_{3}$}--($(L1)+(90:0.5*\s)$)
(R1)node[right]{$a_{1}$}--($(R1)+(90:0.5*\s)$)
(L2)node[right]{$a_{7}$}--($(L2)+(150:0.5*\s)$)
(L2)--($(L2)+(210:0.5*\s)$)
(L3)node[left]{$a_{6}$}--($(L3)+(-120:0.5*\s)$)
(B)node[above]{$a_{9}$}--($(B)+(-120:0.5*\s)$)
(B)--($(B)+(-60:0.5*\s)$)
(R3)node[right]{$a_{8}$}--($(R3)+(-60:0.5*\s)$)
(R2)node[left]{$a_{2}$}--($(R2)+(30:0.5*\s)$)
(R2)--($(R2)+(-30:0.5*\s)$)
(H)node[below right]{$a_{4}$}--($(H)+(90:0.5*\s)$)
;

\draw[line width = 1.5pt]
(L3)--(B)
(R3)--(R2)
(R1)--(H)
(L1)--(L2);

\node[regular polygon, regular polygon sides=3, inner sep = 1, fill, draw] () at (O) {};
\node[rectangle, inner sep = 2.5, fill, draw] () at (H) {};
\node[regular polygon, regular polygon sides=3, inner sep = 1, fill, draw] () at (L1) {};
\node[rectangle, inner sep = 2.5, fill, draw] () at (L2) {};
\node[rectangle, inner sep = 2.5, fill, draw] () at (L3) {};
\node[regular polygon, regular polygon sides=3, inner sep = 1, fill, draw] () at (R1) {};
\node[rectangle, inner sep = 2.5, fill, draw] () at (R2) {};
\node[rectangle, inner sep = 2.5, fill, draw] () at (R3) {};
\node[rectangle, inner sep = 2.5, fill, draw] () at (B) {};
\end{tikzpicture}}
\subcaptionbox{}[0.45\linewidth]{\begin{tikzpicture}[line width = 1pt]
\def\s{1}
\coordinate (O) at (0, 0);
\coordinate (LC) at (144:\s);
\coordinate (H) at ($(LC)!1!72:(O)$);
\coordinate (L1) at ($(LC)!1!72:(H)$);
\coordinate (L2) at ($(LC)!1!72:(L1)$);
\coordinate (L3) at ($(LC)!1!72:(L2)$);
\coordinate (RC) at (36:\s);
\coordinate (R1) at ($(RC)!1!-72:(H)$);
\coordinate (R2) at ($(RC)!1!-72:(R1)$);
\coordinate (R3) at ($(RC)!1!-72:(R2)$);
\coordinate (B) at ($(R3)!1!36:(O)$);

\path [draw=red, postaction={on each segment={mid arrow=red}}]
(O)node[above right]{$a_{5}$}--(H)
(O)--(R3)--(B)
(L3)--(B)
(H)--(L1)
(L2)--(L3)
(R2)--(R1);
\path [draw=red, postaction={on each segment={end arrow=red}}]
(L1)node[left]{$a_{7}$}--($(L1)+(60:0.5*\s)$)
(L1)--($(L1)+(120:0.5*\s)$)
(R1)--($(R1)+(60:0.5*\s)$)
(R1)node[right]{$a_{2}$}--($(R1)+(120:0.5*\s)$)
(L2)node[right]{$a_{3}$}--($(L2)+(180:0.5*\s)$)
(L3)node[left]{$a_{4}$}--($(L3)+(-120:0.5*\s)$)
(B)node[above]{$a_{9}$}--($(B)+(-120:0.5*\s)$)
(B)--($(B)+(-60:0.5*\s)$)
(R3)node[right]{$a_{8}$}--($(R3)+(-60:0.5*\s)$)
(R2)node[left]{$a_{1}$}--($(R2)+(0:0.5*\s)$)
(H)node[below right]{$a_{6}$}--($(H)+(90:0.5*\s)$)
;

\draw[line width = 1.5pt]
(O)--(L3)
(L1)--(L2)
(H)--(R1)
(R2)--(R3);

\node[regular polygon, regular polygon sides=3, inner sep = 1, fill, draw] () at (O) {};
\node[rectangle, inner sep = 2.5, fill, draw] () at (H) {};
\node[rectangle, inner sep = 2.5, fill, draw] () at (L1) {};
\node[regular polygon, regular polygon sides=3, inner sep = 1, fill, draw] () at (L2) {};
\node[rectangle, inner sep = 2.5, fill, draw] () at (L3) {};
\node[rectangle, inner sep = 2.5, fill, draw] () at (R1) {};
\node[regular polygon, regular polygon sides=3, inner sep = 1, fill, draw] () at (R2) {};
\node[rectangle, inner sep = 2.5, fill, draw] () at (R3) {};
\node[rectangle, inner sep = 2.5, fill, draw] () at (B) {};
\end{tikzpicture}}
\subcaptionbox{}[0.45\linewidth]{\begin{tikzpicture}[line width = 1pt]
\def\s{1}
\coordinate (O) at (0, 0);
\coordinate (LC) at (144:\s);
\coordinate (H) at ($(LC)!1!72:(O)$);
\coordinate (L1) at ($(LC)!1!72:(H)$);
\coordinate (L2) at ($(LC)!1!72:(L1)$);
\coordinate (L3) at ($(LC)!1!72:(L2)$);
\coordinate (RC) at (36:\s);
\coordinate (R1) at ($(RC)!1!-72:(H)$);
\coordinate (R2) at ($(RC)!1!-72:(R1)$);
\coordinate (R3) at ($(RC)!1!-72:(R2)$);
\coordinate (B) at ($(R3)!1!36:(O)$);

\path [draw=red, postaction={on each segment={mid arrow=red}}]
(H)--(O)
(O)node[above right]{$a_{5}$}--(L3)
(O)--(R3)--(B)
(L1)--(H)
(L3)--(L2)
(R2)--(R1);
\path [draw=red, postaction={on each segment={end arrow=red}}]
(L1)node[left]{$a_{3}$}--($(L1)+(90:0.5*\s)$)
(R1)node[right]{$a_{2}$}--($(R1)+(60:0.5*\s)$)
(R1)--($(R1)+(120:0.5*\s)$)
(L2)--($(L2)+(150:0.5*\s)$)
(L2)node[right]{$a_{7}$}--($(L2)+(210:0.5*\s)$)
(L3)node[left]{$a_{6}$}--($(L3)+(-120:0.5*\s)$)
(B)node[above]{$a_{9}$}--($(B)+(-120:0.5*\s)$)
(B)--($(B)+(-60:0.5*\s)$)
(R3)node[right]{$a_{8}$}--($(R3)+(-60:0.5*\s)$)
(R2)node[left]{$a_{1}$}--($(R2)+(0:0.5*\s)$)
(H)node[below right]{$a_{4}$}--($(H)+(90:0.5*\s)$)
;

\draw[line width = 1.5pt]
(L3)--(B)
(R3)--(R2)
(R1)--(H)
(L1)--(L2);

\node[regular polygon, regular polygon sides=3, inner sep = 1, fill, draw] () at (O) {};
\node[rectangle, inner sep = 2.5, fill, draw] () at (H) {};
\node[regular polygon, regular polygon sides=3, inner sep = 1, fill, draw] () at (L1) {};
\node[rectangle, inner sep = 2.5, fill, draw] () at (L2) {};
\node[rectangle, inner sep = 2.5, fill, draw] () at (L3) {};
\node[rectangle, inner sep = 2.5, fill, draw] () at (R1) {};
\node[regular polygon, regular polygon sides=3, inner sep = 1, fill, draw] () at (R2) {};
\node[rectangle, inner sep = 2.5, fill, draw] () at (R3) {};
\node[rectangle, inner sep = 2.5, fill, draw] () at (B) {};
\end{tikzpicture}}
\subcaptionbox{}[0.45\linewidth]{\begin{tikzpicture}[line width = 1pt]
\def\s{1}
\coordinate (O) at (0, 0);
\coordinate (LR2) at (0, 2*\s);
\coordinate (LC) at (144:\s);
\coordinate (H) at ($(LC)!1!72:(O)$);
\coordinate (L1) at ($(LC)!1!72:(H)$);
\coordinate (L2) at ($(LC)!1!72:(L1)$);
\coordinate (L3) at ($(LC)!1!72:(L2)$);
\coordinate (RC) at (36:\s);
\coordinate (R1) at ($(RC)!1!-72:(H)$);
\coordinate (R2) at ($(RC)!1!-72:(R1)$);
\coordinate (R3) at ($(RC)!1!-72:(R2)$);
\coordinate (B) at ($(R3)!1!36:(O)$);

\path [draw=red, postaction={on each segment={mid arrow=red}}]
(H)node[below right]{$a_{3}$}--(O)
(O)node[above right]{$a_{4}$}--(L3)
(O)--(R3)--(B)
(L1)node[below]{$a_{2}$}--(H)
(L3).. controls ($(180:0.5*\s)+(L3)$) and ($(180:3*\s)+(LR2)$) ..(LR2)
(R1)--(LR2)
(LR2)node[above]{$a_{6}$}.. controls ($(0:3*\s)+(LR2)$) and ($(0:0.5*\s)+(R3)$) ..(R3);
\path [draw=red, postaction={on each segment={end arrow=red}}]
(L1)--($(L1)+(0:0.5*\s)$)
(R1)node[below]{$a_{1}$}--($(R1)+(180:0.5*\s)$)
(L3)node[left]{$a_{5}$}--($(L3)+(-120:0.5*\s)$)
(B)node[above]{$a_{8}$}--($(B)+(-120:0.5*\s)$)
(B)--($(B)+(-60:0.5*\s)$)
(R3)node[right]{$a_{7}$}--($(R3)+(-60:0.5*\s)$)
(H)--($(H)+(90:0.5*\s)$)
;

\draw[line width = 1.5pt]
(L3)--(B)
(R1)--(H)
(L1)--(LR2);

\node[regular polygon, regular polygon sides=3, inner sep = 1, fill, draw] () at (O) {};
\node[rectangle, inner sep = 2.5, fill, draw] () at (H) {};
\node[regular polygon, regular polygon sides=3, inner sep = 1, fill, draw] () at (L1) {};
\node[rectangle, inner sep = 2.5, fill, draw] () at (LR2) {};
\node[rectangle, inner sep = 2.5, fill, draw] () at (L3) {};
\node[regular polygon, regular polygon sides=3, inner sep = 1, fill, draw] () at (R1) {};
\node[rectangle, inner sep = 2.5, fill, draw] () at (R3) {};
\node[rectangle, inner sep = 2.5, fill, draw] () at (B) {};
\end{tikzpicture}}
\caption{A local matching and orientation.}
\label{XX}
\end{figure}
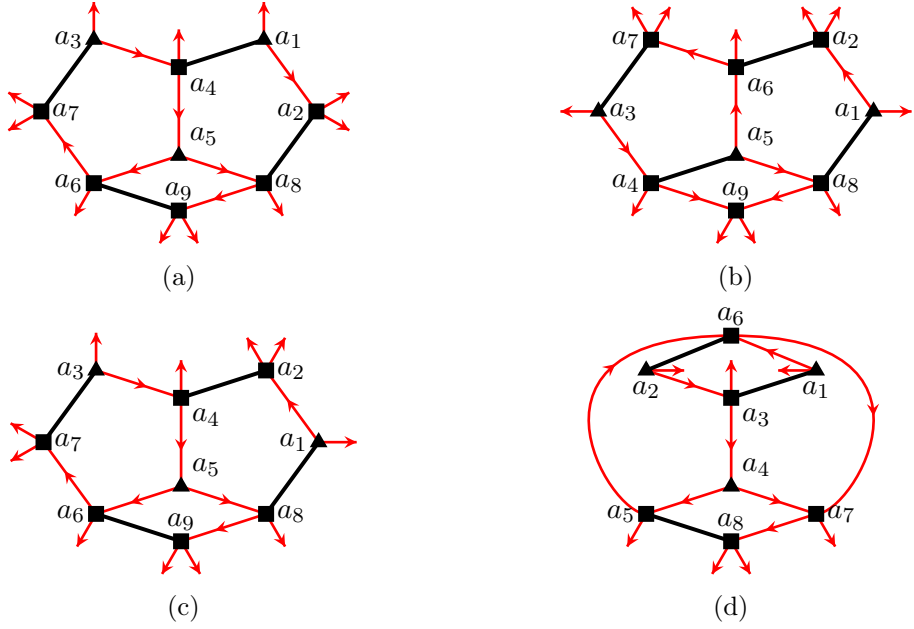

\ref{6ADJACENT3}. Assume $f$ is a $6$-face. If the boundary of $f$ is not a $6$-cycle, then it consists of two triangles, and in that case, no $3$-face is adjacent to $f$. Then we assume the boundary of $f$ is a $6$-cycle $[v_{1}v_{2}\dots v_{6}]$, and a $3$-face $g = [uv_{1}v_{2}]$ is adjacent to $f$. Since $\delta(G) \geq 3$, we have $u \notin \{v_{3}, v_{6}\}$. If $u \in \{v_{4}, v_{5}\}$, then there exists a subgraph isomorphic to \cref{CommonC}, which is a contradiction. Hence, $g$ and $f$ must be normally adjacent, which contradicts the absence in \cref{AT48D}. 

\ref{3444FACE}. Consider a $(3, 4, 3, 4)$-face $[b_{1}b_{2}b_{3}b_{4}]$, and let $X = \{b_{1}, b_{2}, b_{3}, b_{4}\}$. Note that $[b_{1}b_{2}b_{3}b_{4}]$ has no chords. By minimality assumption, $G - X$ can be $(2, 1)$-decomposed into $(D', M')$. Let $M$ be the union of $M'$ and $\{b_{1}b_{2}, b_{3}b_{4}\}$, and let $D$ be obtained from $D' \cup \{\overrightarrow{b_{1}b_{4}}, \overrightarrow{b_{3}b_{2}}\}$ by directing all edges with one endpoint in $X$ and the other endpoint in $\overline{X}$ from $X$ to $\overline{X}$. Then $(D, M)$ is a $(2, 1)$-decomposition of $G$, a contradiction.

\ref{LIGHT3VERTEX}. Let $X$ be the set of vertices on the three faces. By minimality assumption, $G - X$ can be $(2, 1)$-decomposed into $(D', M')$. Now, we label vertices as depicted in \cref{XX}. The bold edges represent a matching $M''$. For each other edge $a_{i}a_{j}$ incident with $X$, we direct it by $\overrightarrow{a_{i}a_{j}}$ whenever $i < j$ or $a_{j} \notin X$. Let $M$ be the union of $M'$ and $M''$, and let $D$ be the resulting orientation of $G - M$. Then $(D, M)$ forms a desired $(2, 1)$-decomposition.

\ref{ReducibleConfiguration}. Let $X$ be the vertex set of any subgraph as depicted in \cref{REDUCIBLEFIGURE}. Then $G - X$ can be $(2, 1)$-decomposed into $(D', M')$. Let $D''$ be the directed edges in \cref{DM}, and $M''$ be the set of bold edges in \cref{DM}. Let $M$ be the union of $M'$ and $M''$ and $D$ be the union of $D'$ and $D''$. Therefore, $(D, M)$ is a desired $(2, 1)$-decomposition of $G$.
\end{proof}

We begin with a primary charge function $\omega$, which is described as $\omega(z) = \deg(z) - 4$ for each $z \in V \cup F$. The Handshaking Theorem and Euler's formula ensure that 
\[
\sum_{z \in V \cup F} \omega(z) = 0.
\]
After the discharging process is finished, we obtain a new charge function $\omega'$. It is important to notice that the total charges can be preserved throughout this process. However, we will prove that $\omega'(z) \geq 0$ for all $z \in F \cup V$, and there exists an element $z \in F \cup V$ with $\omega'(z) > 0$, which leads to a contradiction.

\begin{enumerate}[label=\bf R\arabic*., ref=R\arabic*]
\item\label{R1} Every $3$-face gets $\frac{1}{3}$ from every adjacent face.
\item\label{R2} Consider a $3$-vertex $w$ incident with $h_{1}, h_{2}$ and $h_{3}$. If $h_{1}$ has size at most $4$, then $w$ gets $\frac{1}{2}$ from each of $h_{2}$ and $h_{3}$; otherwise $w$ gets $\frac{1}{3}$ from each of $h_{1}, h_{2}$ and $h_{3}$. 
\item\label{R3} Assume $x$ is a vertex with $\deg(x) \geq 5$. Then $x$ gives $\frac{1}{6}$ to each incident $4^{+}$-face. Furthermore, if $x$ is on a $3$-face $g=[xyz]$ and $yz$ is on another face $g'$, then $x$ gives $\frac{1}{6}$ to $g'$ through $g$. 
\end{enumerate}

\emph{Every vertex concludes with a nonnegative final charge}. Let $u$ be an arbitrary vertex with $\deg(u) = 3$. If $u$ is on a $4^{-}$-face, then the remaining two faces have sizes at least $5$, and $\omega'(u) = \omega(u) + 2 \times \frac{1}{2} = 0$ due to \ref{R2}. If $u$ is not on any $4^{-}$-face, then each incident face is a $5^{+}$-face, which results in $\omega'(u) = \omega(u) + 3 \times \frac{1}{3} = 0$ due to \ref{R2}. 

It is observed that $4$-vertices do not participate in the discharging procedure, thus $\omega'(u) = \omega(u) = 0$. If $\deg(u) \geq 5$, then $\omega'(u) = \omega(u) - \deg(u) \times \frac{1}{6} = \frac{5\deg(u) - 24}{6} > 0$ by \ref{R3}. 

\emph{Every face concludes with a nonnegative final charge}. Let $f = [x_{1}x_{2} \dots x_{d}]$ denote a $d$-face, and let $t = |\{x_{i} \mid \deg(x_{i}) = 3\}|$. Then we can deduce that $t\leq \frac{d}{2}$ from \cref{LEM}\ref{a3}. According to \cref{LEM}\ref{a3}, $f$ is adjacent to at most $d - t$ triangular-faces.

{\bf Case 1. \bm{$d = 3$}.} Then no incident face is a $4^{-}$-face, thus $f$ gets $\frac{1}{3}$ from each adjacent face, resulting in $\omega'(f) = \omega(f) + \frac{1}{3} \times 3 = 0$ due to \ref{R1}. 

{\bf Case 2. \bm{$d = 4$}.} Then $f$ does not send out any charge, implying $\omega'(f) \geq \omega(f) = 0$. 

{\bf Case 3. \bm{$d = 5$}.} By invoking \cref{LEM}\ref{a3}, $f$ contains at most two $3$-vertices. If $f$ is not adjacent to any $3$-face, then $f$ just sends charge to incident $3$-vertices; thus, $\omega'(f) \geq \omega(f) - 2 \times \frac{1}{2} = 0$ due to \ref{R2}. Next, consider the case that $f$ is adjacent to a $3$-face. According to \cref{LEM}\ref{5ADJACENT3} \& \ref{53ADJACENT4}, $f$ is adjacent to precisely one $3$-face $f^{*}$ and zero $4$-faces. If $t \leq 1$, then $f$ sends at most $\frac{1}{2}$ to vertices, and $\frac{1}{3}$ to faces; thus, $\omega'(f) \geq \omega(f) - \frac{1}{2} - \frac{1}{3} > 0$ due to \ref{R1} and \ref{R2}. Now let $t = 2$ and $f^{*} = [ux_{1}x_{2}]$ be the $3$-face. If one of $x_{1}, \dots, x_{5}$ and $u$ is a $5^{+}$-vertex, then $f$ gets at least $\frac{1}{6}$ from these $5^{+}$-vertices, implying $\omega'(f) \geq \omega(f) + \frac{1}{6} - \frac{1}{3} - (\frac{1}{3} + \frac{1}{2}) = 0$. So we may assume that $t = 2$ and no vertex in $\{u, x_{1}, \dots, x_{5}\}$ is a $5^{+}$-vertex. Hence, $f$ contains three $4$-vertices and two $3$-vertices. Suppose that $\deg(x_{1}) = \deg(x_{2}) = 4$. It follows that $\deg(x_{3}) = \deg(x_{5}) = 3$ and $\deg(x_{4}) = 4$. Since the configuration depicted in \cref{53Fig3} is reducible, we have $\deg(u) = 4$. As a result, the configuration shown in \cref{53Fig1} will appear. Assume by symmetry that $x_{2}$ is a $3$-vertex and $u, x_{1}, x_{3}$ are $4$-vertices. If $\deg(x_{5}) = 4$ and $\deg(x_{4}) = 3$, then the configuration as depicted in \cref{53Fig2} occurs, which a contradiction. If $\deg(x_{5}) = 3$ and $\deg(x_{4}) = 4$, then the configuration as depicted in \cref{53Fig3} occurs, which is a contradiction. 

{\bf Case 4. \bm{$d = 6$}.} Then every incident face is a $4^{+}$-face by \cref{LEM}\ref{6ADJACENT3}. Thus, $\omega'(f) \geq \omega(f) - 3 \times \frac{1}{2} > 0$.

{\bf Case 5. \bm{$d \geq 7$}.} Therefore, $\omega'(f) \geq \omega(f) - (d - t) \times \frac{1}{3} - t \times \frac{1}{2} = \frac{2d}{3} - 4 - \frac{t}{6} \geq \frac{7d}{12} - 4 > 0$. 

If $G$ has an element $x$ in $F \cup V$ with $\omega'(x) > 0$, then we complete the proof. So we may assume that $\omega'(x) = 0$ for every element $x \in F \cup V$. By the above arguments, there are no $5^{+}$-vertices and no $6^{+}$-faces. Since there are no adjacent $4^{-}$-faces, there must exist some $5$-faces. However, as we discussed in Case~3, $5$-faces and $3$-faces are not adjacent, and there are two minor $3$-vertices on a $5$-face. Let $f_{1} = [x_{1}x_{2}x_{3}x_{4}x_{5}]$ be a $5$-face adjacent to a $4$-face $g = [x_{2}x_{3}wu_{3}]$ with $\deg(x_{2}) = 3$. By \cref{LEM}\ref{3444FACE}, $f$ must be a $(3, 4, 4, 4)$-face. Since there are no adjacent $4^{-}$-faces, the third face $f_{2}$ incident with $x_{2}$ is a $5$-face. We may assume that $f_{2} = [x_{1}x_{2}u_{3}u_{4}u_{5}]$, and $f_{2}$ is incident with two minor $3$-vertices. However, this contradicts \cref{LEM}\ref{LIGHT3VERTEX}.

\section*{Declarations}
The authors declare that they have no known competing financial interests or personal relationships that could have appeared to influence the work reported in this paper.


\end{document}